\documentclass[12pt]{article}

\usepackage{amsmath,amsthm}
\usepackage{amssymb}
\usepackage{epsfig}

\usepackage[utf8x]{inputenc}
\usepackage{easybmat}
\usepackage{graphicx,hyperref}

\usepackage{url}

\newtheorem{theorem}{Theorem}[section]

\theoremstyle{definition}
\newtheorem{define}[theorem]{Definition}
\newtheorem{remark}[theorem]{Remark}

\begin{document}

\title{A Four-Body Convex Central Configuration\\ 
with Perpendicular Diagonals Is Necessarily a Kite}

\author{Montserrat Corbera\thanks{Departament de Tecnologies Digitals i de la Informaci\'{o}, 
Universitat de Vic, montserrat.corbera@uvic.cat}
\and
Josep M. Cors\thanks{Departament de Matem\`{a}tiques,
Universitat Polit\`{e}cnica de Catalunya, cors@epsem.upc.edu}
\and  Gareth E. Roberts\thanks{
Dept. of Mathematics and Computer Science,
College of the Holy Cross, groberts@holycross.edu}}

\maketitle

\begin{abstract}
We prove that any four-body convex central configuration with perpendicular diagonals must be a kite configuration.  
The result extends to general power-law potential functions, including the planar four-vortex problem.
\end{abstract}



{\bf Key Words:}  Central configuration, $n$-body problem, $n$-vortex problem

\section{Introduction}

Central configurations are an important class of solutions in $n$-body problems.   They lead directly to homographic motions, where the 
initial shape of the configuration is preserved throughout the orbit, and play an important role in the study of the topology of the integral manifolds.
In applied settings, central configurations have proved useful for designing low-cost low-energy space missions~\cite{mars}
and have been discerned in numerical simulations of the eyewall in hurricanes~\cite{davis}.

Locating a central configuration involves solving a challenging set of nonlinear algebraic equations.  
One approach to making the problem more tractable is to impose a geometric
constraint on the shape of the configuration.  In~\cite{pitu-gr}, Roberts and Cors investigated four-body co-circular central configurations,
where the bodies are assumed to lie on a common circle.  Other constraints employed involve symmetry, such as assuming the configuration
has an axis of symmetry~\cite{albouy}, or that it consists of nested regular $n$-gons~\cite{montse}.

Recently,
Li, Deng and Zhang have shown that the diagonals of an isosceles trapezoid central configuration are not perpendicular~\cite{Li}.
We extend this result further by proving that the {\em only} four-body convex central configurations 
with perpendicular diagonals are the kite configurations (symmetric with respect to a diagonal).
Here, ``convex'' means that no body is contained in the convex hull of the other three bodies.  
We prove this result for general power-law potential functions as well as for the planar four-vortex problem.
It is hoped that the techniques described in this paper will prove fruitful for tackling similar
open problems in celestial mechanics~\cite{albouy-pblms}.

\section{Four-Body Planar Central Configurations}
\label{sec:ccc}

We begin by deriving an important equation for four-body central configurations involving mutual distances.
Good references for this material are Schmidt~\cite{schmidt}, Hampton, Roberts, and Santoprete~\cite{HRS}
(for the vortex case), and the recent book chapter by Moeckel~\cite{rick-book}.

Let $q_i \in \mathbb{R}^2$ and $m_i$ denote the position and mass, respectively, of the
$i$-th body.  Except for the case of $n$ point vortices, we will assume that the masses are positive.
Denote $r_{ij} = ||q_i - q_j||$ as the distance between the $i$-th and $j$-th bodies.
If $M = \sum_{i=1}^n m_i$ denotes the sum of the masses, then the {\em center of mass} is given
by $c = \frac{1}{M} \sum_{i=1}^n m_i q_i$.  The motion of the bodies is governed by the 
potential function 
$$
U_\alpha(q) \; = \;  \sum_{i < j}^n  \; \frac{m_i m_j}{r_{ij}^{\alpha}},
$$
where $\alpha > 0$ is a parameter.  The classical $n$-body problem corresponds to $\alpha = 1$.
The {\em moment of inertia} with respect to the center of mass, which measures the relative size of the configuration, is given by
\begin{equation}
I(q) \; = \;  \sum_{i=1}^{n} \; m_i \| q_i -  c \|^2 \; = \;
\frac{1}{M} \sum_{i < j}  \;  m_i m_j r_{ij}^2 .
\label{eq:inertia}
\end{equation}

There are many equivalent ways of describing a central configuration.  We follow the topological 
approach.  Let $I_0$ be some positive constant.  

\begin{define}
A planar {\em central configuration} $(q_1, \ldots, q_n) \in \mathbb{R}^{2n}$
is a critical point of $U_\alpha$ subject to the constraint $I = I_0$.
\end{define}

It is important to note that, due to the invariance of $U_\alpha$ and $I$ under isometries, any rotation, translation, or
scaling of a central configuration still results in a central configuration.

\begin{figure}[htb]
 \centering
 \includegraphics[height=7.5cm,keepaspectratio=true]{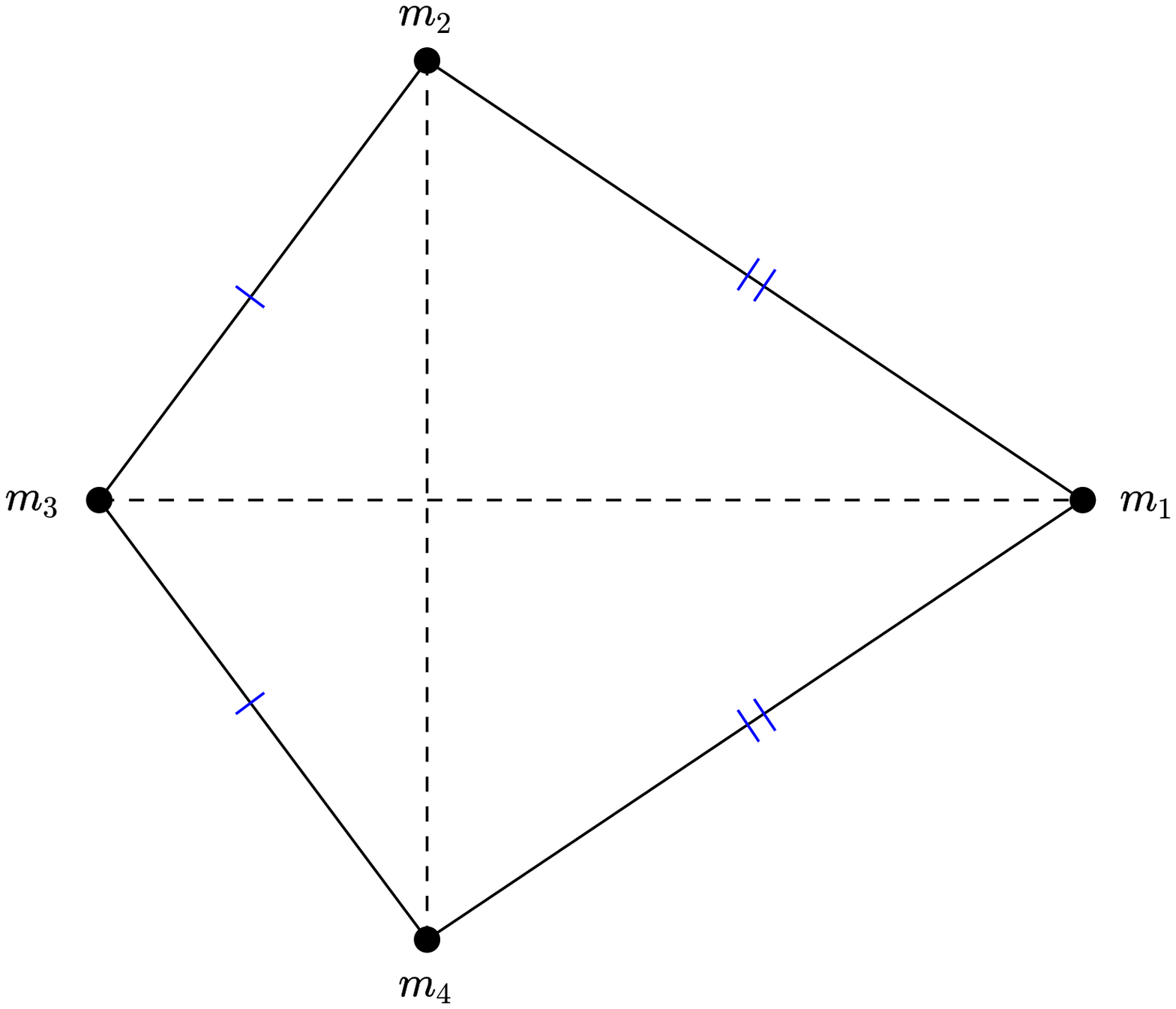}
 \hspace{0.05in}
 \includegraphics[height=8.4cm,keepaspectratio=true]{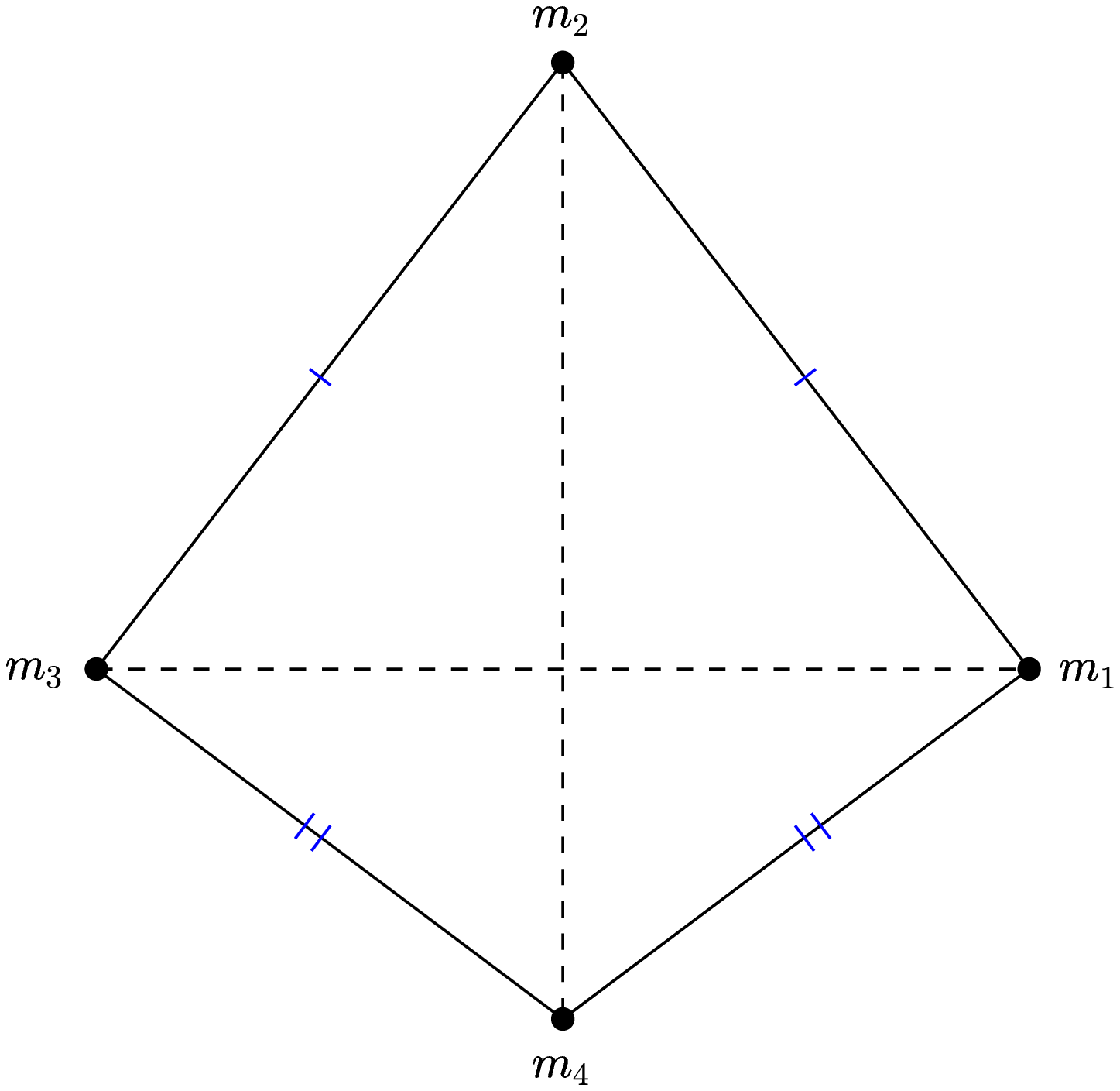}
 \caption{Two convex kite configurations with different symmetry axes.}
 \label{Fig:kites} 
\end{figure}

We now restrict to central configurations in the planar four-body problem.   A 
configuration is {\em convex} if no body lies inside or on the convex hull of the other three bodies
(e.g., a rhombus or a trapezoid); otherwise, it is {\em concave}.  
A {\em kite} configuration contains two bodies on an axis of symmetry and two
bodies located symmetrically with respect to this axis (see Figure~\ref{Fig:kites}). The diagonals of a convex kite configuration are perpendicular.  
The main result of this paper is to prove that kite configurations are the only central configurations with this property.

\begin{theorem}
If the diagonals of a four-body convex central configuration are perpendicular, then the configuration must be a kite.
\label{mainThm}
\end{theorem}

For four bodies, the six mutual distances $r_{12}, r_{13}, r_{14}, r_{23}, r_{24}, r_{34}$ are excellent
coordinates.  These distances describe an actual configuration in the plane if the
Cayley-Menger determinant, defined as
$$
V \; = \;  \left|  \begin{array}{ccccc}
		0 & 1  & 1   & 1  & 1 \\
	         1 & 0  & r_{12}^2 & r_{13}^2 & r_{14}^2  \\
                  1 & r_{12}^2 &  0  & r_{23}^2  & r_{24}^2 \\
                  1 & r_{13}^2 & r_{23}^2 & 0 & r_{34}^2    \\  
                  1 & r_{14}^2 & r_{24}^2 & r_{34}^2 & 0 
              \end{array} \right|,
$$
vanishes and the triangle inequality $r_{ij} + r_{jk} > r_{ik}$ holds for any choice of
indices with $i \neq j \neq k$.   We must impose the constraint $V=0$ to find planar central configurations; 
otherwise, the only critical points of $U_\alpha$ restricted
to $I = I_0$ are regular tetrahedra.  Thus, we seek critical points of the function
\begin{equation}
U_\alpha + \lambda(I - I_0) + \mu V
\label{eq:cps}
\end{equation}
satisfying $I = I_0$ and $V = 0$, where $\lambda$ and $\mu$ are Lagrange multipliers.

Let $A_i$ be the oriented area of the triangle whose vertices contain all bodies
except for the $i$-th body.  Assuming the bodies in a convex quadrilateral are ordered sequentially, we have
$A_1, A_3 > 0$ and $A_2, A_4 < 0$.  An important formula involving the Cayley-Menger determinant is
\begin{equation}
\frac{ \partial V}{ \partial r_{ij}^2} \; = \;  - 32 \, A_i  A_j \, .
\label{eq:cayley-deriv}
\end{equation}
Formula~(\ref{eq:cayley-deriv}) is only valid when restricting to planar configurations.

Differentiating function~(\ref{eq:cps}) with respect to $r_{ij}$ and applying formula~(\ref{eq:cayley-deriv}) yields 
\begin{equation}
m_i m_j (s_{ij} - \lambda^{'}) \; = \;  \sigma A_i A_j,
\label{eq:cc}
\end{equation}
where $s_{ij} = r_{ij}^{-(\alpha + 2)}, \lambda^{'} = 2\lambda/(\alpha M),$ and $\sigma = -64 \mu/\alpha$.
Group the six equations defined by~(\ref{eq:cc}) as follows:
\begin{equation}
\begin{split}
& m_1 m_2(s_{12} - \lambda') = \sigma A_1 A_2, \qquad   m_3  m_4(s_{34} - \lambda') = \sigma A_3 A_4,\\
& m_1 m_3(s_{13} - \lambda') = \sigma A_1 A_3, \qquad   m_2  m_4(s_{24} - \lambda') = \sigma A_2 A_4,\\
& m_1 m_4(s_{14} - \lambda') = \sigma A_1 A_4, \qquad   m_2  m_3(s_{23} - \lambda') = \sigma A_2 A_3.
\end{split}
\label{eq:ccGroup}
\end{equation}
Multiplying the equations together pairwise and cancelling the common terms yields the well-known Dziobek relation
\begin{equation}
(s_{12} - \lambda^{'})(s_{34} - \lambda^{'}) \; = \;  (s_{13} - \lambda^{'})(s_{24} - \lambda^{'})\;  = \;  (s_{14} - \lambda^{'})(s_{23} - \lambda^{'}).
\label{equs:dzio}
\end{equation}
Finally, after eliminating $\lambda^{'}$ from equation~(\ref{equs:dzio}), we obtain the important relation
\begin{equation}
(r_{24}^\beta - r_{14}^\beta)(r_{13}^\beta - r_{12}^\beta) (r_{23}^\beta - r_{34}^\beta)  \; = \;  (r_{12}^\beta - r_{14}^\beta) (r_{24}^\beta - r_{34}^\beta) (r_{13}^\beta - r_{23}^\beta),
\label{Eq:consist}
\end{equation}
where $\beta = \alpha + 2$.

Equation~(\ref{Eq:consist}) is a necessary condition for a four-body planar central configuration.  
Assuming that $V = 0$, it is also sufficient, although further restrictions are needed to ensure that the masses
are positive.   For a convex configuration, the two diagonals must each be strictly longer than any of the exterior
sides.  For instance, if the bodies are ordered sequentially, then we have $r_{13}, r_{24} > r_{12}, r_{14}, r_{23}, r_{34}$.
Moreover, the shortest and longest exterior sides must lie opposite each other (with equality only in the case
of the square).  Thus, if $r_{12}$ is the length of the longest exterior side, then $r_{12} \geq r_{14}, r_{23} \geq r_{34}$.
These geometric relations can be derived by dividing different pairs of equations in~(\ref{eq:ccGroup})
and using the convexity of the configuration (see~\cite{schmidt} for details).

\section{Proof of Theorem~\ref{mainThm}}

\begin{figure}[h!]
 \centering
 \includegraphics[height=9cm,keepaspectratio=true]{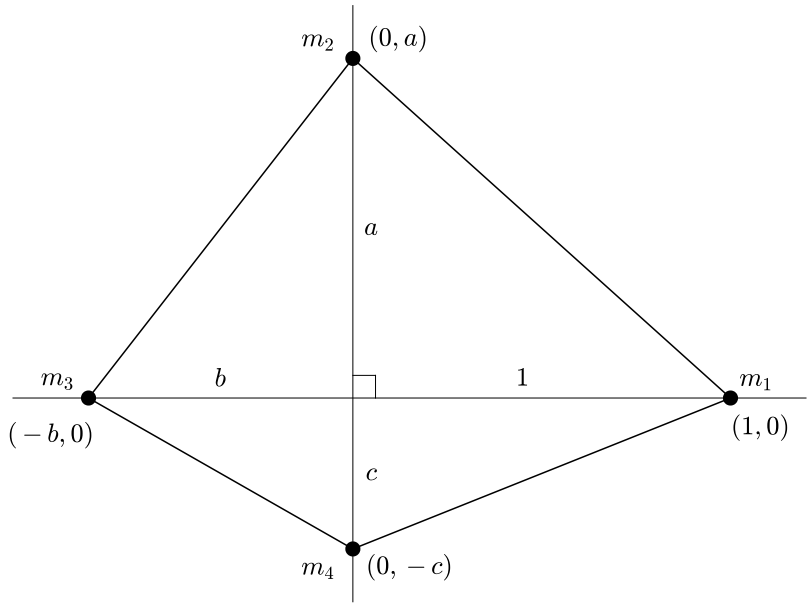}
 \vspace{-0.2in}
 \caption{Coordinates for the proof of Theorem~\ref{mainThm}.}
 \label{Fig:setup} 
\end{figure}

Without loss of generality, we can assume that the bodies in the convex central configuration are ordered sequentially 
in a counterclockwise fashion around the configuration.  Thus, the lengths of the diagonals are $r_{13}$ and $r_{24}$, and these are
larger than any of the four exterior sides $r_{12}, r_{23}, r_{14}, r_{34}$.  We can also assume that the bodies are labeled so
that $r_{12}$ is the longest exterior side.    Since the diagonals are perpendicular, we can apply a rotation and translation 
to place the bodies on the coordinate axes.  Finally, by scaling the configuration, we can assume that $q_1 = (1,0)$.  
Let the other positions be given by $q_2 = (0, a), q_3 = (-b, 0)$ and $q_4 = (0, -c)$, where $a, b, c$ are positive real variables (see Figure~\ref{Fig:setup}).
In these coordinates kite configurations occur in two distinct cases:  $a = c$ (horizontal axis of symmetry; left plot in Figure~\ref{Fig:kites}), or
$b = 1$ (vertical axis of symmetry; right plot in Figure~\ref{Fig:kites}).   If both $a = c$ and $b = 1$, the configuration is a rhombus.

\vspace{-0.2in}

The mutual distances for the exterior sides are given by
$$
r_{12}^2 = a^2 + 1,  \quad  r_{23}^2 = a^2 + b^2,  \quad  r_{34}^2 = b^2 + c^2,  \quad  \mbox{ and } \quad r_{14}^2 = 1 + c^2,
$$
while the lengths of the two diagonals are simply $r_{13} = 1 + b$ and $r_{24} = a + c$.  Since $r_{12}$ is the length of the longest exterior
side, we have $r_{12} \geq r_{14}$, which implies $a \geq c$ and $r_{23} \geq r_{34}$.

Denote $F$ as the function
$$
F(a,b,c) \; = \;  (r_{24}^\beta - r_{14}^\beta)(r_{13}^\beta - r_{12}^\beta)(r_{23}^\beta - r_{34}^\beta) -  (r_{12}^\beta - r_{14}^\beta)(r_{24}^\beta - r_{34}^\beta)(r_{13}^\beta - r_{23}^\beta),
$$
where each mutual distance is treated as a function of the variables $a, b, c$.  Let $\Gamma$ be the level surface $F(a,b,c) = 0$
restricted to the octant where $a$, $b$, and $c$ are strictly positive and satisfy the inequalities
\begin{equation}
r_{13}, r_{24} \;  >  \;  r_{12} \; \geq \; r_{14}, r_{23} \; \geq r_{34}.
\label{geomCond}
\end{equation}
By equation~(\ref{Eq:consist}), $\Gamma$ is equivalent to the set of convex central configurations with perpendicular diagonals, positive masses, and our
choice of labeling.

It is easy to see that the convex kite configurations are contained in $\Gamma$.  
For instance, if $b = 1$, then $r_{12} = r_{23}$ and $r_{14} = r_{34}$, from which it
follows that $F(a,1,c) = 0$ for any choice of $a$ or $c$.  Likewise, if $a = c$, then $r_{12} = r_{14}, r_{23} = r_{34}$, and $F(a,b,a) = 0$
for any choice of $a$ or $b$.  We will show that these are the {\em only} two solutions to the equation $F = 0$ on $\Gamma$.

Using equation~(\ref{Eq:consist}) and the partial derivatives $\partial r_{13}/\partial b = 1,  \partial r_{23}/\partial b = b/r_{23}$, 
and $\partial r_{34}/\partial b = b/r_{34}$, we compute that
\begin{eqnarray*}
\left. \frac{\partial F}{\partial b} \right|_\Gamma   & = &  (r_{24}^\beta - r_{14}^\beta) [\beta r_{13}^{\beta - 1} (r_{23}^\beta - r_{34}^\beta) + b\beta(r_{13}^\beta - r_{12}^\beta)(r_{23}^\alpha - r_{34}^\alpha)]   \\
  &  &  - (r_{12}^\beta - r_{14}^\beta)[-b\beta r_{34}^\alpha (r_{13}^\beta - r_{23}^\beta) + \beta(r_{24}^\beta - r_{34}^\beta)(r_{13}^{\beta - 1} - b r_{23}^\alpha)]  \\[0.13in]
  & = &  \beta r_{13}^{\beta - 1}[(r_{24}^\beta - r_{14}^\beta)(r_{23}^\beta - r_{34}^\beta) - (r_{12}^\beta - r_{14}^\beta)(r_{24}^\beta - r_{34}^\beta)] \\[0.05in]
  &  &  + \, b \beta(r_{24}^\beta - r_{14}^\beta) (r_{13}^\beta - r_{12}^\beta)(r_{23}^\alpha - r_{34}^\alpha) \\[0.05in]
  &  &  + \, b \beta(r_{12}^\beta - r_{14}^\beta)[ r_{34}^\alpha(r_{13}^\beta - r_{23}^\beta) + r_{23}^\alpha(r_{24}^\beta - r_{34}^\beta)] \\[0.13in]
  & = &  \frac{\beta r_{13}^{\beta - 1} (r_{12}^\beta - r_{14}^\beta)(r_{24}^\beta - r_{34}^\beta)(r_{12}^\beta - r_{23}^\beta)}{r_{13}^\beta - r_{12}^\beta}  
      + b \beta(r_{24}^\beta - r_{14}^\beta) (r_{13}^\beta - r_{12}^\beta)(r_{23}^\alpha - r_{34}^\alpha) \\[0.05in]
     &  &  + \, b \beta (r_{12}^\beta - r_{14}^\beta)[ r_{34}^\alpha (r_{13}^\beta - r_{23}^\beta) + r_{23}^\alpha (r_{24}^\beta - r_{34}^\beta)] .
\end{eqnarray*}

Recall that $a \geq c$ on $\Gamma$.  Fix $a$ and $c$ such that $a > c$ and look for a value of $b$ that satisfies $F(a,b,c) = 0$.  
One such value is $b = 1$, a kite configuration.  But $a > c$ implies that $r_{23} - r_{34} > 0$, 
from which it follows that $\partial F/\partial b > 0$ on $\Gamma$ due to the inequalities in~(\ref{geomCond}).
Since $F$ is a differentiable function, it follows that there can be only one solution to $F = 0$, and this must be $b = 1$.
Note that if $a = c$, then $r_{12} = r_{14}$ and $r_{23} = r_{34}$, which yields $\partial F/\partial b = 0$, as expected (the plane
$a = c$ satisfies $F = 0$ for any value of~$b$).  Thus, $b = 1$ or $a = c$ are the only solutions to $F = 0$ on $\Gamma$.
This completes the proof.

\begin{remark}
The surface $\Gamma$ (kite configurations with our particular choice of labeling), consists of a subset of the union of two orthogonal planes $b=1$ and $a=c$.
In order to satisfy the inequalities in~(\ref{geomCond}), on the plane $b=1$ we have $1/\sqrt{3} < a < \sqrt{3}$ and
$\sqrt{a^2 + 1} - a < c \leq a$, while on the plane $a=c$ we have $1/\sqrt{3} < a < \sqrt{3}$ and
$\sqrt{a^2 + 1} - 1 < b \leq 1$.  The two planes intersect in a line that corresponds to the one-parameter family of rhombii central configurations.
It is straight-forward to check that 0 is a critical value for $F(a,b,c)$, as the partial derivatives all vanish 
at points of the form $(a,1,a)$.  This is consistent with the fact that $\Gamma$, defined by the pre-image $F^{-1}(0)$,
is not a manifold.  
\end{remark}

\begin{remark}
For completeness, we use the equations in~(\ref{eq:ccGroup}) to compute the 
masses for each type of kite configuration.  Without loss of generality, we can scale the masses so that $m_1 = 1$.
For kite configurations with $b = 1$ (symmetric with respect to the $y$-axis), we have $m_3 = m_1 = 1$, 
$$
m_2 \; = \;  \frac{2c}{a+c} \cdot \frac{s_{14} - s_{13}}{s_{23} - s_{24}}, \quad \mbox{ and } \quad
m_4 \; = \;  \frac{2a}{a+c} \cdot \frac{s_{23} - s_{13}}{s_{14} - s_{24}}.
$$
For kite configurations with $a = c$ (symmetric about the $x$-axis), we have $m_1 = 1$, 
$$
m_2 \; = \;  m_4 \; = \;  \frac{b+1}{2b} \cdot \frac{s_{14} - s_{13}}{s_{23} - s_{24}}, \quad \mbox{ and } \quad
m_3 \; = \;  \frac{1}{b} \cdot \frac{(s_{14} - s_{13}) (s_{14} - s_{24})}{ (s_{23} - s_{13})(s_{23} - s_{24})}.
$$
\end{remark}

\section{The Vortex Problem}

Next we consider the case of point vortices with arbitrary circulations $\Gamma_i \in \mathbb{R}$.  Unlike the $n$-body problem,
it is possible that some of the circulations are negative.   Central configurations in the planar $n$-vortex problem are
critical points of the Hamiltonian $H = -\sum_{i<j} \Gamma_i \Gamma_j \ln (r_{ij})$ subject to the constraint $I = I_0$,
where $I$ is given by equation~(\ref{eq:inertia}), but is now called the {\em angular impulse}.  See Chapter~2 of~\cite{newton}
for a good overview of the planar $n$-vortex problem.

As before, in order to locate planar four-body central configurations,
we must impose the additional constraint $V = 0$ using the Cayley-Menger determinant.
The argument used earlier to derive equation~(\ref{Eq:consist}) works equally well 
in this setting.  The mutual distances must satisfy  
\begin{equation}
(r_{24}^2 - r_{14}^2)(r_{13}^2 - r_{12}^2) (r_{23}^2 - r_{34}^2)  \; = \; (r_{12}^2 - r_{14}^2) (r_{24}^2 - 
r_{34}^2) (r_{13}^2 - r_{23}^2),
\label{Eq:consistVorts}
\end{equation}
although there are no inequalities such as~(\ref{geomCond}) restricting the mutual distances
since the circulations can be of opposite sign.

Suppose that we have a convex central configuration of four vortices with perpendicular diagonals.  
Using the same setup and variables as in the proof of Theorem~\ref{mainThm}, 
equation~(\ref{Eq:consistVorts}) becomes
$$
(a^2 + 2ac - 1)(b^2 + 2b - a^2)(a^2 - c^2) - (a^2 - c^2)(a^2 + 2ac - b^2)(1 + 2b - a^2) \; = \;  0,
$$
which factors nicely into
$$
2(a^2 - c^2)(b^2 - 1)(ac + b) \; = \;  0.
$$
Since $a, b, c$ are strictly positive, we quickly deduce that $a = c$ or $b = 1$, and the configuration 
must be a kite.  We have proven the following theorem, the analog to Theorem~\ref{mainThm} for the 
point vortex problem.

\begin{theorem}
If the diagonals of a four-vortex convex central configuration are perpendicular, then the configuration must be a kite.
\end{theorem}

\begin{remark}
As with the vortex case, choosing $\alpha$ to be an even natural number leads to nice factorizations of equation~(\ref{Eq:consist}).
For example, if $\alpha = 2$ (the strong force potential), then equation~(\ref{Eq:consist}) factors as 
$$
(a^2 - c^2)(b^2 - 1) \cdot p_2(a,b,c) \; = \; 0,
$$
where $p_2$ is an eighth-degree polynomial with 47 terms whose coefficients are all positive.  Similarly, for the case $\alpha = 4$, 
equation~(\ref{Eq:consist}) factors as 
$$
(a^2 - c^2)(b^2 - 1) \cdot p_4(a,b,c) \; = \; 0,
$$
where $p_4$ is a fourteenth-degree polynomial with 210 terms whose coefficients are all positive.  These factorizations were quickly obtained using Maple~\cite{maple}.
Since $a,b,c$ are strictly positive, $a=c$ or $b=1$ follows immediately, and the argument used in the proof of Theorem~\ref{mainThm}
is not needed.
\end{remark}

\vspace{0.2in}

\noindent  {\bf Acknowledgments:} 
G. Roberts was supported by a grant from the National Science Foundation (DMS-1211675).
M. Corbera and J. M. Cors were partially supported by MINECO grant MTM2013-40998-P;
J. M. Cors was also supported by AGAUR grant 2014 SGR 568.

\bibliographystyle{amsplain}

\end{document}